\documentclass[12pt]{amsart}
\usepackage{enumerate}
\usepackage{mathrsfs,amsthm,amsfonts,amsmath,amssymb,amsrefs}
\usepackage{mathtools}
\mathtoolsset{showonlyrefs,showmanualtags}

\usepackage{color}
\usepackage{tikz}
\usepackage[colorlinks=true, pdfstartview=FitV, linkcolor=blue, citecolor=blue, urlcolor=blue]{hyperref}
\usepackage[top=1.5in, bottom=1.5in, left=1.25in, right=1.25in]	{geometry}
\usepackage{url}

\newtheorem{theorem}{Theorem}[section]

\newtheorem{definition}[theorem]{Definition}
\newtheorem{lemma}[theorem]{Lemma}
\newtheorem{proposition}[theorem]{Proposition}
\newtheorem{conjecture}[theorem]{Conjecture}

\def\la{\langle}
\def\ra{\rangle}

\def\one{\mathbf 1}

\numberwithin{equation}{section}

\begin{document}

\title{Weighted Estimates for One Sided Martingale Transforms}

\author[Wei Chen]{Wei Chen}
\address{ School of Mathematical Sciences, Yangzhou University, Yangzhou 225002, China; School of Mathematics, Georgia Institute of Technology, Atlanta GA 30332, USA}
\email {weichen@yzu.edu.cn}

\author{Rui Han}
\address{ School of Mathematics, Georgia Institute of Technology, Atlanta GA 30332, USA}
\email {rui.han@math.gatech.edu}

\author[M. T. Lacey] {Michael T. Lacey}   

\address{ School of Mathematics, Georgia Institute of Technology, Atlanta GA 30332, USA}
\email {lacey@math.gatech.edu}

\thanks{
Wei Chen is supported by the National Natural Science Foundation of China(11771379), the Natural Science Foundation of Jiangsu Province(BK20161326), and the Jiangsu Government Scholarship for Overseas Studies(JS-2017-228). 
The research of R. Han is supported in part by grant from the US National Science Foundation, DMS-1800689. 
The research of M. T. Lacey is supported in part by grant from the US National Science Foundation, DMS-1600693, and the
Australian Research Council ARC DP160100153.}

\begin{abstract}
Let $ Tf =\sum_{  I} \varepsilon_I \langle f,h_{I^+}\rangle h_{I^-}$. Here,  $ \lvert  \varepsilon _I\rvert=1 $,
and $ h_J$ is the Haar function defined on dyadic interval $ J$.  We show that, for instance,
\begin{equation*}
\lVert T \rVert _{L ^{2} (w) \to L ^{2} (w)}  \lesssim
[w]  _{A_2 ^{+}}  .
\end{equation*}
Above, we use the one sided $ A_2$ characteristic for the weight $ w$.
This is an instance of a one sided $A_2$ conjecture.
Our proof of this fact is difficult, as the very quick known proofs of the $A_2$ theorem do
not seem to apply in the one sided setting.
\end{abstract}


\maketitle
\section{One-sided martingale transform}
We prove new sharp one sided weight inequalities for certain kinds of one sided martingale transforms.
One sided weights are variants of the usual $ A_p$ weights on the real line, for which the following supremum is finite:
\begin{equation}\label{e:Ap}
[w] _{A_p ^{+}} = \sup _{I}  \frac{w (I^-)} {\lvert  I^-\rvert } \Bigl[ \frac{\sigma (I^+)} {\lvert  I^+\rvert } \Bigr] ^{p-1},
\end{equation}
where $ w$ is a non-negative function,   locally integrable, and $ \sigma = w ^{ - \frac{1} {p-1}}$
is also locally integrable.  The set $ I = (a, a+ 2 \delta )$ is an interval, with left and right halves  $ I ^{-} = (a, a+ \delta )$,
and $ I ^{+} = (a+ \delta , a+ 2 \delta )$ respectively.
There are variants of this notion in higher dimension which we will not directly discuss.

For the endpoint cases of \( p= 1, \infty \), we use the one sided maximal function
\begin{equation*}
M_+ f = \sup _{I } \frac {\mathbf 1_{I^-} } {\lvert  I^+\rvert } \int _{I^+} \lvert  f \rvert\; dx .
\end{equation*}
Define \( M _{-}\) similarly. We set
\begin{equation}\label{e:1}
[w] _{A_1^+} = \Bigl\lVert  \frac { M_- w}{ w} \Bigr\rVert_ \infty ,
\qquad
[w] _{A_ \infty ^{+}} = \sup _{I} w (I) ^{-1} \int _{I} M_- (w\mathbf 1_{I})\; dx.
\end{equation}
The definitions of $ A _{p} ^{-}$ follow similarly, and it is important to note the duality between these expressions.  In particular,
$
[w ] _{A_p ^{+}} = [ \sigma ] _{A_{p'} ^{-}} ^{p-1}
$.

The qualitative aspects of the $ A_p^+$ theory closely matches that of the usual $ A_p$ theory.
And, in the latter, the study of sharp constants, which has been under very rapid development in the $ A_p$ theory.  But some of the corresponding results for $ A_p ^{+}$ seem  much harder to establish.
The maximal function estimates were established by Sawyer \cite{MR849466}, and we have this result, matching known results for the usual maximal function.

\begin{theorem}\label{t:Max}  \cite{MR3338703}*{Thm1.5 and Thm 1.9}  These inequalities hold:
\begin{align}\label{e:MaxWeak}
\lVert M_+ \rVert _{L ^{p} (w) \to L ^{p, \infty } (w)} & \lesssim [w]^{1/p} _{A^+ _{p}}  , \qquad 1< p < \infty ,
\\ \label{e:MaxStrong}
\lVert M_+ \rVert _{L ^{p} (w) \to L ^{p} (w)} & \lesssim  ([w] _{A^+ _{p}} [\sigma ] _{A_  \infty  ^-} ) ^{1/p}  , \qquad 1< p < \infty .
\end{align}
\end{theorem}

We remark that we have
\begin{equation}\label{e:sigma}
[ \sigma ] _{A_ \infty  ^{-}} \leq [\sigma ] _{A_{p'} ^{-}} = [w] _{A_p ^{+}} ^{p'-1}.
\end{equation}

The main contribution of this paper is to establish some sharp $ A _p ^{+}$ estimates for
(maximal truncations of) one sided martingale transforms.  The latter are defined by
\begin{align*}
Tf &=\sum_{\textup{dyadic}\, I} \varepsilon_I \langle f,h_{I^+}\rangle h_{I^-},
\\
T _{\sharp }f &= \sup _{\delta  >0}
\Bigl\lvert \sum_{\substack{\textup{dyadic}\, I\\ \lvert  I\rvert > \delta    }} \varepsilon_I \langle f,h_{I^+}\rangle h_{I^-} \Bigr\rvert,
\end{align*}
where
$\varepsilon_I$ takes values in $\{-1, 1\}$ and
$h$ is the Haar function defined by
$$h_I(x)=\frac {\one_{I^+}(x)-\one_{I^-}(x)} {\sqrt{|I|}}.$$
Our main result is the sharp  weak $ L ^{p}$ inequality, for all $ 1< p < \infty $,
and the sharp inequality on $ L ^2 $.

\begin{theorem}\label{thm:main}  The following inequalities are uniform over all  one sided martingale transforms.
\begin{align}\label{e:Tweak}
\lVert T _{\sharp } \rVert _{L ^{p} (w) \to L ^{p, \infty } (w)} & \lesssim [w]^{\frac{1}{p}} _{A_p ^{+}} [w ] _{A_ \infty ^{+}} ^{1/p'}, \qquad 1< p < \infty ,
\\  \label{e:A2T}
\lVert T \rVert _{L ^{2} (w) \to L ^{2} (w)} & \lesssim
{[w]^{\frac{1}{2}} _{A_2 ^{+}}
\max \bigl\{ [\sigma ]  _{A_ \infty ^{-}} , [w] _{A_ \infty^{+}}\}}^{\frac{1}{2}}.
\end{align}
\end{theorem}

This is a very modest assertation, as compared to  the advanced state of the $ A_p$ theory.
We will address what we think it true, and some of the unexpected complications we encountered,
in the concluding section of this paper.

Concerning the proof, it is noteworthy that the proof of the maximal function estimates \eqref{e:MaxStrong}
are relatively simple.  One uses the weak-type result, together with a sharp reverse H\"older estimate
 \cite{MR3338703}*{Thm 1.8}.
Marcinkiewicz interpolation finishes the proof.  (This type of argument was first identified by Buckley \cite{MR1124164}.)

No such argument can work for the martingale transforms, since the $ L ^{\infty }$ endpoint estimate is too large.
Indeed, most of the strategies that are so successful for the usual $ A_p$ weights do not seem to generalize to the one sided setting.
We adapt methods from
\cites{MR2993026,MR2657437} in order to complete the proof.
Namley, we establish a general version of a distributional inequality, which is the key to plus/minus characteristic of the $ A_p$ condition.  This is combined with testing conditions, and a corona construction.
There are new complications in carrying out this program,
and some serious obstructions to proving what one would expect to be true.

\bigskip
The theory of $ A_p ^+$ weights was started by
Sawyer \cite{MR849466}. He showed the natural analog of Muckenhoupt's famous theorem.
The $ A_p ^{+}$ condition is characterized by those positive \textup{a.e.} weights $ w$ for which $ M_+$
is bounded on $ L ^{p} (w)$.
Amir, Foranzi and Mart\'{i}n-Reyes \cite{MR1376747} established several natural analogs of
the $ A_p ^+$ theory for singular integrals.
The recent paper of Chill and Kr\'{o}l \cite{MR3811173} gives an elegant development  of the $ A_p ^{+}$ theory,
especially with an eye towards applications in parabolic PDEs.  While there is a sizable literature on
$ A_p ^+$ weights, even in the singular integral setting, the results we could find were qualitiative in nature.
As far as we are aware, the question of sharp constants has not been  addressed before for
singular integral like objects in the $ A_p ^+$ setting.

\section{Background} 

It is a useful remark that the one weight inequality \( \lVert T f\rVert _{L ^{p} (w)} \lesssim \lVert f\rVert _{L ^{p} (w)}\) is equivalent to the two weight inequality
\begin{equation} \label{e:2wt}
\lVert T(\sigma f)\rVert _{L ^{p} (w)} \lesssim \lVert f\rVert _{L ^{p} (\sigma )}, \qquad
\sigma = w ^{1- p'}.
\end{equation}
The latter inequality is a two weight inequality, and is convenient as it dualizes correctly,
namely the inequality above is equivalent to  $ \lVert T^\ast (w \phi )\rVert _{L ^{p'} (\sigma )} \lesssim \lVert \phi \rVert _{L ^{p'} (w )}$.

There is a characterizations of the two weight inequality \eqref{e:2wt} in the case of $ p=2$, due to
Nazarov, Treil and Volberg \cite{MR2407233}. It completely solves the question for dyadic operators.
For the sake of clarity, we recall it, in the context of our one sided martingale transforms.

\begin{theorem}\label{t:NTV}   \cite{MR2407233} Let $ (w, \sigma )$ be a pair of weights, and $ T $ a one sided martingale transform.  We have $  \lVert T (\sigma \cdot )\rVert _{L ^2 (\sigma ) \to L ^2 (w)} \lesssim \mathfrak  T$,
where $ \mathfrak T$ is the best constant in the inequalities below, uniformly over intervals $ I$.
\begin{align}
\lVert \mathbf 1_{I} T( \sigma \mathbf 1_{I})\rVert _{L ^2 (w)} &\leq \mathfrak  T \sigma (I) ^{1/2},
\\   \label{e:test}
\lVert \mathbf 1_{I} T ^{\ast} ( w \mathbf 1_{I})\rVert _{L ^2 (\sigma )} &\leq \mathfrak  T w (I) ^{1/2}.
\end{align}
\end{theorem}

The condition \eqref{e:test} is referred to as a \emph{testing condition.}
It is a fundamental reduction in complexity in proving the $ L ^2 $ bound.

A corresponding theorem in the $ L ^{p}$ setting is necessarily more complicated.
A true characterization has not yet been found, but must be of a vector valued nature, see \cite{MR3612889}.
Indeed, the main result of this last paper is powerful, but  are stated with an assumption about `quadratic $ A_p$ condition',  \cite{MR3612889}*{(3.2)}, which is ill suited to the one sided case.

We do not need a full characterization, and if one restricts attention to dyadic Calder\'on-Zygmund operators,  very sharp conditions can be given.
We will use this corollary to   \cite{MR2993026}*{Thm 4.3}. It is a two weight inequality for
maximal truncations of one sided martingale transforms.

\begin{theorem}\label{t:twoWeight}  Let \( w, \sigma \) be two weights, and \( T\) a one sided martingale transform.  We have the weak-type bound
\begin{equation*}
\lVert T _{\sharp } (\sigma f)\rVert _{L ^{p,\infty} (w)} \lesssim (\mathfrak M^+_p + \mathfrak T_p) \lVert f\rVert _{L ^{p} (\sigma )}
\end{equation*}
where the two constants on the right are the best constants in the following inequalities, holding
for all functions \( f\), and intervals \( I\).
\begin{align}\label{e:Mtest}
\lVert M _{+} (\sigma f) \rVert _{L ^{p,\infty} (w)} &\leq \mathfrak M^+_p \lVert f\rVert _{L ^{p} (\sigma )},
\\ \label{e:Ttest}
\int_{I} T _{\sharp } (\sigma \mathbf 1_{I} f) \; dw  &\leq \mathfrak T_p  \lVert f\rVert _{L ^{p} (\sigma )} w(I) ^{1/p'}.
\end{align}
\end{theorem}

The Theorem as formulated above is not a corollary, since the cited \cite{MR2993026}*{Thm 4.3}
has the usual maximal function appearing in \eqref{e:Mtest}.
The proof however extends immediately to the version above.
We will comment more of this in the final section of the paper.

The condition \eqref{e:Ttest} does not look at all like the corresponding testing condition \eqref{e:test}.  Let us explain here
why they are similar.
The maximal truncation operator is not linear.
But its boundedness is equivalent to the boundedness of a family of linear operators given
as follows.
For any measurable $ \delta \;:\; \mathbb R \to (0, \infty )$, let us set
\begin{equation*}
T _{\delta } f (x)= \sum_{\substack{I \;:\; \lvert  I\rvert > \delta (x) \\ }}
\langle f , h _{I^+} \rangle h _{I ^{-}} (x)
\end{equation*}
be a linearization of the maximal truncations.
It is easy to see that $ T _{\sharp}$ is bounded if and only if the family of linear operators $ T _{\delta }$
satisfy a norm bound independent of the choice of measurable $ \delta $.

Using linearity, note that the integral in \eqref{e:Ttest}  can be written as
\begin{equation*}
\langle    T _{\delta }(f \mathbf 1_{I} \sigma ),   \varphi _I \cdot w  \rangle \leq
 \mathfrak T_p  \lVert f\rVert _{L ^{p} (\sigma )} w(I) ^{1/p'}.
\end{equation*}
This inequality should hold uniformly over all choices of truncation $ \delta $, and $
 \lvert  \varphi _I\rvert = \mathbf 1_{I} $.  Therefore, the condition \eqref{e:Ttest} is the same as
\begin{equation}\label{e:TT}
\lVert \mathbf 1_{I} T _{\delta } ^{\ast} (w \varphi _I ) \rVert _{L ^{p'} (\sigma )} \leq
 \mathfrak T_p   w(I) ^{1/p'}.
\end{equation}
This is the form of \eqref{e:Ttest} that we will use.  It is a testing condition.

The appearance of $ T _{\delta } ^{\ast} $ is not so familiar.  Crucially, these operators  satisfy a weak $ L ^{1} $ inequality.
We state here, and remark that the proof is not easy.

\begin{proposition}\label{p:weak} \cite {MR2993026}*{Thm. 9.3}
$  T _{\delta } ^{\ast}$ maps $ L ^{1}$ to weak $ L ^{1}$, uniformly over the choice of the one sided martingale transform and the choice of the linearization, $\delta $.

\end{proposition}

Finally, we need a Lemma which states that a weak type distributional estimate is sufficient for a John-Nirenberg type estimate.

\begin{lemma} (\cite{MR3176607}*{Lemma 5.5}, \cite {MR2993026}*{Lemma 10.2})
\label{l:JN}
Let $ \{ \phi _{I } \,:\, I\in \mathcal D\} $ be a collection of functions indexed by dyadic intervals which
are  supported on $ I$, and
constant on the  grandchildren of $I$.
Suppose that there is a constant $C$ so that
for any dyadic interval $I_{0}$, and any collection $\mathcal E$ of
dyadic subintervals $I\subset I_{0}$, there holds
\begin{equation}   \label{e:weak}
\Bigl\lvert \Bigl \{   \Bigl\lvert \sum _{I\in \mathcal E} \phi_{I } \Bigr\rvert > C \Bigr\} \Bigr\rvert
< \tfrac 12 \lvert I_{0} \rvert .
\end{equation}
Then, we have
\begin{equation}  \label{e:JN}
\Bigl\lvert \Bigl \{   \Bigl\lvert \sum _{I\in \mathcal E} \phi_{I } \Bigr\rvert > (C+1) \lambda  \Bigr\} \Bigr\rvert
<  2 ^{(1-\lambda)/2} \lvert I_{0} \rvert , \qquad  \lambda > 1.
\end{equation}
\end{lemma}

\section{The Distributional Lemma} 

The key component is a one sided variant of a distributional lemma discovered in \cite{MR2657437}.
Fix an interval $ I_0$, and for integers $ a \in \mathbb Z $,
let $ \mathcal K _{a}$ be those intervals
dyadic $I\subseteq I_0^-$ for which we have
\begin{gather}\label{e:a}
2 ^{a} <  \langle \sigma  \rangle _{I ^{+}} ^{p-1} \langle w \rangle _{I ^{-}} \leq 2 ^{a+1},
\\  \label{e:w<}
\langle w \rangle _{I ^{-}} \leq 2 \langle w \rangle _{I_0 ^{-}}.
\end{gather}

\begin{lemma}\label{l:Distribution}  For an absolute constant $c>0$ and constant $C_p>0$ that depends (only) on $p$, we have
\begin{equation}\label{e:Distribution}
\sigma  \Bigl(  \bigl\lvert T _{\delta  , \mathcal K_a} ^{\ast}  (w \varphi  _{I_0})\bigr\rvert > C_p \lambda \langle w  \rangle _{I_0 ^{-}} \Bigr) \lesssim  
\begin{cases}
\lambda^{-\frac{2p'}{p+1}}  \cdot \sigma  (I_0 ^{+}),\ \ \text{for } 0<\lambda<1,\\
\\
e^{-c\lambda}\cdot \sigma(I_0^+),\ \ \text{for } \lambda\geq 1.
\end{cases}
\end{equation}
Above, for any collection $\mathcal K$,
\begin{equation}  \label{e:K*}
(T _{\delta  , \mathcal K} ^{\ast}(f))(x):=
\sum _{\substack{I\in \mathcal K}} \varepsilon _I
\langle f, h _{I^-}  \mathbf 1_{ \lvert  I\rvert > \delta (x) } \rangle  h _{I^+}(x).
\end{equation}
The function $ \varphi _ I$ satisfies $ \lvert  \varphi _I\rvert = \mathbf 1_{I} $.
\end{lemma}

The point of this inequality is that it holds in the two weight setting, provided the $ A_p ^{+}$ product is approximately constant. 
Moreover, one take care to note that the set on the left is in $ I _0 ^{-}$, and its $\sigma$-measure is controlled by $\sigma(I_0^+)$.

\begin{proof}
The collection $ \mathcal K _{a}$ is further divided into collections $ \mathcal K _{a} ^{b}$,
for $ b\in \mathbb N $, where  $ I\in \mathcal K _{a} ^{b}$ if $ I\in \mathcal K_{a}$ and in addition,
\begin{equation}
\label{e:b}
2 ^{-b} \langle w  \rangle _{I_0 ^{-}} < \langle w  \rangle _{I ^{-}} \leq 2 ^{1-b} \langle w  \rangle _{I_0 ^{-}}.
\end{equation}
By \eqref{e:w<}, this is a decomposition of $ \mathcal K_{a}$.

We show that for constants $c, C>0$, and $\lambda>0$,
\begin{equation}\label{e:D}
\sigma  \Bigl(  \bigl\lvert T _{\delta  , \mathcal K_a^b} ^{\ast}  (w \varphi  _{I_0}) \bigr\rvert >  C\lambda 2 ^{1-b} \langle w  \rangle _{I_0 ^{-}} \Bigr) \lesssim  2^{b /(p-1) - c \lambda } \sigma (I_0 ^{+}).
\end{equation}
Nevertheless, this is a strong condition, and another easy sub additivity argument
completes the proof of \eqref{e:Distribution}.  

Indeed, for $0<\lambda<1$. Since $ \sum_{b=0} ^{\infty } 2 ^{-b(1-\frac{p+1}{2p})} =1/(1-2^{\frac{1-p}{2p}})=:C_p<\infty$, we have
\begin{align*}
\sigma  \Bigl(  \lvert T _{\delta  , \mathcal K_a} ^{\ast}  (w \varphi  _{I_0})\rvert  >  2C C_p \lambda \langle  w \rangle _{I_0 ^{-}} \Bigr)
\leq 
&\sum_{b=0} ^{\infty } \sigma  \Bigl(  {T _{\delta , \mathcal K_a^b} ^{\ast}} (w \varphi  _{I_0}) >   C\lambda 2 ^{1-b(1-\frac{p+1}{2p})} \langle w  \rangle _{I_0 ^{-}} \Bigr)
\\
\lesssim 
&\sigma  (I_0 ^{+}) \sum_{b=0} ^{\infty } 2^{b /(p-1) - c\lambda 2 ^{\frac{b(p+1)}{2p}}} \lesssim \lambda^{-\frac{2p'}{p+1}} \cdot \sigma  (I_0 ^{+}).
\end{align*}

For $\lambda\geq 1$. Since $\sum_{b=0}^{\infty} 2^{-b/2}<10$, we have
\begin{align*}
\sigma  \Bigl(  \lvert T _{\delta  , \mathcal K_a} ^{\ast}  (w \varphi  _{I_0})\rvert  >  20 C \lambda \langle  w \rangle _{I_0 ^{-}} \Bigr)
\leq 
&\sum_{b=0} ^{\infty } \sigma  \Bigl(  {T _{\delta , \mathcal K_a^b} ^{\ast}} (w \varphi  _{I_0}) >   C\lambda 2 ^{1-b/2} \langle w  \rangle _{I_0 ^{-}} \Bigr)
\\
\lesssim 
&e^{-c\lambda}\cdot \sigma  (I_0 ^{+}) \sum_{b=0} ^{\infty } 2^{b /(p-1) - c 2 ^{b/2}} \lesssim e^{-c\lambda}\cdot \sigma  (I_0 ^{+}).
\end{align*}
This proves \eqref{e:Distribution}.

In general, we know little more than the local integrability condition on $ w \varphi  _{I_0} $.
But we have this consequence of the weak type inequality Proposition \ref{p:weak}:
There is a constant $C$ so that for any subinterval $J\subset I_{0}$ and any collection $\mathcal E \subset \mathcal K_{a} ^{b}$
of subintervals of $J$, we have
\begin{equation} \label{e:JJNN}
\bigl\lvert \bigl\{  \lvert  T _{\delta  , \mathcal E} ^{\ast} ( w \varphi_{J} ) \rvert > C 2^{1-b} \la w \ra_{I_0^-}  \bigr\} \bigr\rvert \leq \tfrac 12 \lvert J\rvert .
\end{equation}
This follows from the dyadic structure, and  condition \eqref{e:b}, which controls the local $ L ^{1} $
norms of $ w \varphi _{I_0}$.

Turn the sufficient conditions for a John-Nirenberg estimate in Lemma~\ref{l:JN}.
The inequality \eqref{e:JJNN} is the hypothesis \eqref{e:weak}.
The other hypotheses of Lemma \ref{l:JN} are easy to check.
It follows that for $\lambda>0$,
\begin{align}\label{eq:JN_KabS}
\left|\left\lbrace x:~  | T _{\delta  , \mathcal K_a^b} ^{\ast} (x) ( w \varphi_{I_{0}} ) |> C \lambda 2^{1-b} \la w\ra_{I_0^-} \right\rbrace \right|
\lesssim 2^{-c\lambda} |I_0|.
\end{align}
Note that for $0<\lambda<1$, this estimate is trivial.

Our task is to convert this Lebesgue measure estimate into $ \sigma $-measure.  By combining conditions \eqref{e:a} and \eqref{e:b}, we have for any $ I \in \mathcal K _{a,b}$
\begin{align*}
\langle \sigma  \rangle _{I ^{+}}  ^{p-1}& \sim  2 ^{a}  \langle w  \rangle _{I ^{-}} ^{-1}
\\
& \sim  2 ^{a + b } \langle w  \rangle _{I ^{-}_0} ^{-1} \sim 2 ^{b}   \langle \sigma   \rangle _{I ^{+}_0} ^{p-1}.
\end{align*}
That is, $\sigma  (I ^{+}) \sim 2^{b/(p-1)} \rho \lvert  I \rvert$, 
with constant $ \rho $ \emph{independent of $ I\in \mathcal K _a^b$}. 

The set on the left in \eqref{eq:JN_KabS} is a union of children of intervals $ I ^{-}$, for $ I\in \mathcal K_a^b$.  
Therefore, \eqref{e:D} follows.  And, the Lemma is proved.
\end{proof}

\section{Proof of Theorem \ref{thm:main}}
There are two assertations, the weak $L^p$ bound, and the $L^2$ bound. Both are proved by
appeal to the two weight theorems, Theorem~\ref{t:twoWeight} and Theorem~\ref{t:NTV}.
By inspection, it suffices to prove the weak $L^p$ bound for the maximal truncations.
This inequality specialized to the case of $ p=2$, and taking duality into account, gives the full $ L ^2 $ result.

From Theorem~\ref{t:twoWeight}, we need to estimate the quantities $ \mathfrak M^+ _{p}$ and $ \mathfrak T _p$ defined in \eqref{e:Mtest} and \eqref{e:Ttest}.  The first $ \mathfrak M^+_p$ is the weak type norm for the maximal function $ M _{+}$. That is given in \eqref{e:MaxWeak}, and is smaller than what is claimed. So we turn to the second constant,
$ \mathfrak T_p$ which is the testing constant for the maximal truncations of  one sided martingale transform.
Namely, the task is to show that uniformly over all dyadic intervals $ I$,
that the inequality \eqref{e:TT} holds.   More exactly, we need to see that
\begin{equation}\label{e:TTT}
\lVert \mathbf 1_{I} T _{\delta } ^{\ast} (w \varphi _I ) \rVert _{L ^{p'} (\sigma )} \leq
 [w ] ^{1/p} _{A_ {p} ^{+}}\, [w]^{{1/p'}} _{A _{\infty } ^{+}}  \cdot  w(I) ^{1/p'}.
\end{equation}

The martingale transform is a sum over all dyadic intervals. We can restrict the sum to those dyadic intervals $ I$ that intersect $ I_0$.  Those $ I$ that  strictly contain $ I_0$ can also be dismissed, since they contribute nothing to the left side of \eqref{e:TT}. The case of $ I=I_0$ is trivial.  Thus, the difficult
case is $ I\subsetneq I_0$, and we restrict attention to this case, by assuming that the multiplier coefficients $ \epsilon _I =0$ unless $ I\subsetneq I_0$.

Our principal tools are the Distribution Lemma~\ref{l:Distribution}  and   a corona type decomposition.
The latter is needed to get to a point that we can apply Lemma~\ref{l:Distribution}.
For an integer $ a$ with $ 2 ^{a} < [w] _{A_p ^{+}}$, let
\begin{align*}
\mathcal K_a &= \{I \subsetneq I_0 \;:\;
2 ^{a} <  \langle \sigma  \rangle _{I ^{+}} ^{p-1} \langle w \rangle _{I ^{-}} \leq 2 ^{a+1}\}.
\end{align*}

Assuming that $ \epsilon _I =0$ unless $ I\in \mathcal K_a$, we will show that
\begin{equation}\label{e:Ta}
\int _{I_0} T _{\delta  , \mathcal K_a} ^{\ast}  (w \varphi  _{I_0})  ^{p'}\; \sigma  (dx)
\lesssim 2 ^{a (p'-1)  } [w ] _{A_ \infty ^{+}} w (I_0).
\end{equation}
Above, we are using the notation of \eqref{e:K*}.
Summing over $ a $ will prove \eqref{e:TTT}.

\begin{definition}[The Corona]\label{def:Cat}
Let
$$C_{a,1}:=\{I\in \mathcal{K}_a:~ I^- \text{ is maximal w.r.t. inclusion} \}.$$
For $I\in \mathcal{C}_{a,1},$ let $\mathcal{C}_{a,2}(I)$ be the ``bad'' children of $I$, defined as below
$$\mathcal{C}_{a,2}(I):=\{J\in \mathcal{K}_a: ~J^-\subset I^-,~\langle w \rangle_{J_-}>2 \langle w \rangle_{I^-}, ~ J^- \text{ is maximal}\},$$
and let
$$\mathcal{C}_{a,2}=\bigcup_{I\in \mathcal{C}_{a,1}}\mathcal{C}_{a,2}(I).$$
We also define $\mathcal{C}_{a,t}$ for $t\geq 3$ inductively.
Let
$$\mathcal{C}_a=\bigcup_{t=1}^{\infty}\mathcal{C}_{a,t}.$$
\end{definition}



We refine the corona decomposition.
\begin{definition}\label{def:minimal_stopping_time}
For $J\in \mathcal{K}_a$, let $J^s$ stand for the minimal stopping interval $S\in \mathcal{C}_a$ with $J^-
\subseteq S^-$.
\end{definition}
The collections below form a partition of $\mathcal{K}_a$.
$$\mathcal{K}_a(S):=\{J\in \mathcal{K}_a:~ J^s=S\},\ \ S\in \mathcal{C}_a,$$


We are now at the core of the argument.  Set $  \tau _S = T _{\delta  , \mathcal K_a (S)} ^{\ast}  (w \varphi  _{I_0})$, and
$$  X_{S,n}=\mathbf 1_{\{ 2 ^{n-1} \la w\ra_{S^-}<|\tau _S|\leq  2 ^{n}\la w\ra_{S^-}\}}\tau _S , \quad n\in \mathbb Z.
$$
Then, estimate
\begin{align}
\bigl\lVert  \mathbf 1_{I_0}  {T^* _{\delta , \mathcal K _{a} }(w \varphi  _{I_0}) }\bigr\rVert _{L ^{p'} (\sigma )}
& =  \Bigl\lVert  \sum_{S\in \mathcal C_a}  \tau _{S}\Bigr\rVert _{L ^{p'} (\sigma )}
\\
& \leq  \sum_{n= - \infty  } ^{\infty }  \Bigl\lVert  \sum_{S\in \mathcal C_a}  X _{S,n} \Bigr\rVert _{L ^{p'} (\sigma )}
\\  \label{e:X}
& \lesssim   \sum_{n= - \infty  } ^{\infty }
\biggl[ \sum_{S\in \mathcal C_a} \lVert  X _{S,n} \rVert _{L ^{p'} (\sigma )} ^{p'} \biggr] ^{1/p'} .
\end{align}
The last inequality follows from the construction of the corona: the values $ \langle w  \rangle _{S ^{-}}$  form a geometric sequence of reals.

Concerning this norms in \eqref{e:X}, the case of $ n \geq n_p$ and $ n <n_p$ are different, here $n_p$ is the least positive integer such that $2^{n-1}\geq C_p$.
For $ n\geq n_p$, we have, again by the construction of the corona and the distribution estimate \eqref{e:Distribution} with $\lambda=2^{n-1}/C_p\geq 1$,
\begin{align*}
\sum_{S\in \mathcal C_a}  \lVert  X _{S,n}\rVert _{L ^{p'} (\sigma )} ^{p'}
& \lesssim  2^{np'}e ^{-c{2^{n-1}}}\sum_{S\in \mathcal C_a}  \langle  w\rangle  _ {S ^{-}} ^{p'} \sigma  (S ^{+})
\\
& \lesssim  2^{np'}e ^{-c{2^{n-1}}} 2 ^{a (p'-1)}\sum_{S\in \mathcal C_a}   \langle  w  \rangle  _ {S ^{-}} \lvert  S ^{+}\rvert
\\
& \lesssim  2^{np'}e ^{-c{2^{n-1}}} 2 ^{a (p'-1)} [w ] _{A ^{+} _{\infty }} w (I_0).
\end{align*}
This estimate decreases rapids in $n$, so that  it is strong enough to imply \eqref{e:Ta}.

The case of $n < n_p$, is very similar. We use the distribution estimate \eqref{e:Distribution} with $\lambda=2^{n-1}/C_p<1$,
\begin{align*}
\sum_{S\in \mathcal C_a}  \lVert  X _{S,n}\rVert _{L ^{p'} (\sigma )} ^{p'}
& \lesssim  2^{np'} 2^{-n\frac{2p'}{p+1}} \sum_{S\in \mathcal C_a}  \langle  w\rangle  _ {S ^{-}} ^{p'} \sigma  (S ^{+})
\\
& \lesssim  2^{np/(p+1)} 2 ^{a (p'-1)}\sum_{S\in \mathcal C_a}   \langle  w  \rangle  _ {S ^{-}} \lvert  S ^{+}\rvert
\\
& \lesssim  2^{np/(p+1)} 2 ^{a (p'-1)} [w ] _{A ^{+} _{\infty }} w (I_0).
\end{align*}
Since $ n \leq n_p$, this is again strong enough to complete the proof of  {\eqref{e:Ta}.
And the proof is finished.\qed

\section{Complements} 
\newcounter{para}
\newcommand\mypara{\par\refstepcounter{para}\textbf{\thepara.}\space}

\mypara It would be very natural to seek a proof of our main theorem using ideas related to
sparse bounds. The latter subject was started in \cite{MR3625108}, and has been very successful.
But we could not find such a bound in the current setting.
That is why we returned to an earlier proof \cite{MR2657437} of the $ A_2$ bound, one that was explored
in \cites{MR2993026,MR2912709,MR3176607}.

\mypara
The source of many of our difficulties seem to be linked to this point.
A key element of the recent developments in the $ A_p$ theory is the following fact:  For any weight $ \mu $
on $ \mathbb R $, the maximal function
\begin{equation*}
M _{\mu } f =   \sup _{I \;:\; \mu (I)>0} \mathbf 1_{I}  \mu (I) ^{-1} \int _{I} f \; d \mu
\end{equation*}
is weakly bounded on $ L ^{1} (\mu )$.
This does not seem to be true for the plus versions of this maximal function.

\mypara The absence of this `universal maximal function' reminds us of a similar obstruction in the multiparameter setting.  Recent results of Barron and Pipher \cite{170905009B} have shown that in fact sparse bounds do not hold in that setting.  Would some version of these arguments hold in the one sided case?

\mypara We pose the question:  What is the simplest natural condition that one can place on a weight $ \mu $ so that the maximal operator
\begin{equation*}
M _{\mu } ^{+}f =   \sup _{I \;:\; \mu (I ^{+})>0} \mathbf 1_{I ^{-}}  \mu (I ^{+}) ^{-1} \int _{I ^{+}} f \; d \mu
\end{equation*}
is weakly bounded on $ L ^{1} (\mu )$?  Or on some $ L ^{p} (\mu )$, for $ 1\leq p < \infty $. Is it $ \mu \in A _{\infty } ^{+}$?

\mypara  Our main theorem admits a straight forward extension to paraproducts, that is operators of the form
\begin{equation*}
T f = \sum_{I} \tau _I  \langle f \rangle _{I ^{+}} h _{I ^{-}},
\end{equation*}
where $ \tau _I$ is a Carleson sequence.
It likewise admits an extension to the setting where the martingale transforms are replaced with
Haar shifts with complexity, in the sense of \cite{MR2912709,MR3204859}.
In the latter case, of complexity, one wants bounds that are polynomial in complexity.
There are versions of our two weight theorems Theorem~\ref{t:NTV} and Theorem~\ref{t:twoWeight}
that account for complexity.
In the interest of clarity, we have not pursued these points here.

\mypara
One of the ways that the classical $ A_p$ theory and the $ A_p ^{+}$  are similar is in the
area of extrapolation.  A robust theory holds in both places.
A quantified version of  the extrapolation of $A_p ^+ $ is described in
Carro, Lorente and Mart\'{i}n-Reyes \cite{MR3771482}.
The extrapolation of strong type norms in the $ A_p^ +$ setting is described in
\cite{MR2797562}.  That discussion is not quantitative, but there should not be a major obstacle to devising such a theory.
With it, one could deduce new strong type inequalities from our inequality, \eqref{e:A2T}.

\mypara
Given the role of two weight inequalities Theorem~\ref{t:NTV} and Theorem~\ref{t:twoWeight}
in our argument, it might be reasonable to investigate these theorems in the setting of one sided
operators.  We comment that \cite{MR2993026}* {Thm 4.7} has a general strong type two weight inequality.
But, like for Theorem~\ref{t:twoWeight}, the discussion is not geared towards one sided operators.
Moreover, the proof of that theorem is not so easy, and the sufficient conditions are substantially more involved.  So, we felt that appealing to the one sided version of that result was not in the spirit of this paper.

\mypara  The paper of Vuorinen \cite{MR3612889} aims for a characterization of the
strong type inequality for general dyadic operators. But, note that the main theorem of this paper
begins with an assumption of a `quadratic  two weight $ \mathcal A _{p}$' condition,
specified in section 3 of that paper.  This condition is not of a one sided nature.  So again, the
Theorem might be attractive to apply in the one sided setting, but it is not completely staight forward to
do so.

\mypara
The overarching conjecture here concerns \emph{one sided Calder\'on-Zygmund operators},
and a one sided version of the main results of \cites{MR2912709,MR3092729}.
These
are operators  $ T $, bounded on $ L ^2 $,  that have a kernel representation
\begin{equation*}
\langle T f, g \rangle = \int K (x,y) f (y) g (x) \; dx\,dy
\end{equation*}
where $ K (x,y)$ satisfies the standard size and smoothness conditions
\begin{align*}
 \lvert  K (x,y)\rvert  & \lesssim  \lvert  x-y\rvert ^{-1} , \qquad x\neq y
 \\
\lvert \nabla ^{\alpha } K (x,y) \rvert  &\lesssim  \lvert  x-y\rvert ^{-1- \alpha },
\end{align*}
where $ 0< \alpha \leq 1$ is fixed. But, one imposes the one sided condition:  $ K (x,y ) = 0 $ if $ x < y$.

\begin{conjecture}\label{j:A2}
Let $ T$ be a one sided Calder\'on-Zygmund operator.  Then, there holds
\begin{equation*}
\lVert T \rVert _{L ^{p} (w)\to L ^{p} (w)}
\lesssim  [w] _{A_p ^{+}} \max \{   [\sigma ] _{A_ \infty  ^{-}}  ^{1/p} ,
[w] _{A _{\infty } ^{+}} ^{1/p'}\} \lesssim [w] _{A_p} ^{\max \{1,  (p-1) ^{-1} \}} .
\end{equation*}
\end{conjecture}

\mypara  One of the beautiful aspects of  the proof of the $ A_2$ theorem in \cite{MR2912709} is the \emph{Hyt\"onen Representation Theorem}, which gives a representation of a Calder\'on-Zygmund operator as a rapidly convergent in complexity
sum of dyadic shifts.  We could not prove the anlagous result in the one sided setting.
\emph{Does  a one sided Hyt\"onen Representation hold?}   We see no reason why it should not hold, but could not find a proof.

\bibliographystyle{alpha,amsplain}	


\end{document}